\documentclass[11pt]{amsart}
\usepackage{amsfonts}
\usepackage{amsmath}
\usepackage{amssymb}
\usepackage{color}
\usepackage{enumerate}
\usepackage{tikz}

\usetikzlibrary{arrows,calc,positioning,trees,shadows}

\newtheorem{theorem}{Theorem}

\newtheorem{corollary}[theorem]{Corollary}

\newtheorem{definition}[theorem]{Definition}

\newtheorem{lemma}[theorem]{Lemma}

\newtheorem{proposition}[theorem]{Proposition}

\def\DB{\overline{\Delta}}

\def\id{id}

\def\NM{{\mathbf N}}

\newenvironment{arb}{\begin{tikzpicture}[baseline,scale=0.5,level distance=7mm,level 1/.style={sibling distance=10mm},level 2/.style={sibling distance=5mm},level 3/.style={sibling distance=3mm},grow=down, font=\scriptsize]
\tikzstyle{ve}=[draw,circle,inner sep=1pt,fill] 
\tikzstyle{vv}=[draw,circle,inner sep=1pt] 
\tikzstyle{vee}=[minimum size=0pt ,inner sep=0pt]}{\end{tikzpicture}}
\newenvironment{arbb}{\begin{tikzpicture}[baseline,scale=0.5,level distance=7mm,level 1/.style={sibling distance=25mm},level 2/.style={sibling distance=5mm},level 3/.style={sibling distance=3mm},grow=down, font=\scriptsize]
\tikzstyle{ve}=[draw,circle,inner sep=1pt,fill] 
\tikzstyle{vv}=[draw,circle,inner sep=1pt] 
\tikzstyle{vee}=[minimum size=0pt ,inner sep=0pt]}{\end{tikzpicture}}

\newcommand{\rd}[1]{\node[ve,label=above:$#1$] {}}

\newcommand{\vl}[1]{node[ve,label=left:$#1$] {}}
\newcommand{\vr}[1]{node[ve,label=right:$#1$] {}}
\newcommand{\vb}[1]{node[ve,label=below:$#1$] {}}

\def\corollaa{\begin{arb}
\rd{(i;i_0)}
child{\vb{i_1}} child[dotted] child[dotted] child{\vb{i_k}};
\end{arb}}

\def\corollab{\begin{arb}
\rd{(i;i_0)}
child{\vb{i_1}}  child{\vb{i_2}};
\end{arb}}

\def\corollac{\begin{arb}
\rd{(i;i_0)}
child{\vl{(i_1 ; i_{1,0})} child{\vb{i_{1,1}}} child[dotted] child[dotted] child{\vb{i_{1,k}}}} child{\vr{i_2}};
\end{arb}}

\def\corollad{\begin{arb}
\rd{(i;i_0)}
child{\vl{i_1}} child{\vr{(i_2 ; i_{2,0})}  child{\vb{i_{2,1}}} child[dotted] child[dotted] child{\vb{i_{2,l}}}} ;
\end{arb}}

\def\corollae{\begin{arbb}
\rd{(i;i_0)}
child{\vl{(i_1 ; i_{1,0})} child{\vb{i_{1,1}}} child[dotted] child[dotted] child{\vb{i_{1,k}}}}child{\vr{(i_2 ; i_{2,0})}  child{\vb{i_{2,1}}} child[dotted] child[dotted] child{\vb{i_{2,l}}}} ;
\end{arbb}}

\def\exarb{\begin{arbb}
\rd{(i;i_0)}
child{\vl{(i_1 ; i_{1,0})} child{\vb{i_{1,1}}} child[missing] child{\vb{i_{1,2}}}} child{\vr{(i_2 ; i_{2,0})}  child{\vb{i_{2,1}}}child[missing] child{\vb{i_{2,2}}}child[missing] child{\vb{i_{2,3}}}} ;
\end{arbb}}

\def\tfdbA{\begin{arb}
\rd{(3;2)}
child{\vb{1}};
\end{arb}}
 
\def\tfdbB{\begin{arb}
\rd{(3;1)}
child{\vb{1}} child{\vb{1}}  ;
\end{arb}}

\def\tfdbC{\begin{arb}
\rd{(3;1)}
child{\vb{2}};
\end{arb}}

\def\tfdbD{\begin{arb}
\rd{(2;1)}
child{\vb{1}};
\end{arb}}

\def\tfdbE{\begin{arb}
\rd{(3;1)}
child{\vr{(2;1)} child{\vb{1}}};
\end{arb}}

\begin{document}

\title[Forest formula]{Right-handed Hopf algebras and the PreLie forest formula.} 

\date{}




\author{Fr\'ed\'eric Menous}
\address{Laboratoire de Math\'ematiques, B\^at. 425,
UMR 8628, CNRS,
Universit\'e Paris-Sud,
91405 Orsay Cedex,
France }
\email{Frederic.Menous@math.u-psud.fr}

\author{Fr\'ed\'eric Patras}
\address{Laboratoire J.A. Dieudonn\'e,
         		UMR 7531, CNRS,
         		Universit\'e de Nice,
         		Parc Valrose,
         		06108 Nice Cedex 2, France.}
\email{Frederic.PATRAS@unice.fr}

\begin{abstract}
Three equivalent methods allow to compute the antipode of the Hopf algebras of Feynman diagrams in perturbative quantum field theory (QFT): the Dyson-Salam formula, the Bogoliubov formula, and the Zimmermann forest formula. Whereas the first two hold generally for arbitrary connected graded Hopf algebras, the third one requires extra structure properties of the underlying Hopf algebra but has the nice property to reduce drastically the number of terms in the expression of the antipode (it is optimal in that sense).

The present article is concerned with the forest formula: we show that it generalizes to arbitrary right-handed polynomial Hopf algebras. These Hopf algebras are dual to the enveloping algebras of preLie algebras -a structure common to many combinatorial Hopf algebras which is carried in particular by the Hopf algebras of Feynman diagrams.
\end{abstract}

\maketitle

\section*{Introduction}
\label{sect:intro}
Three equivalent methods allow to compute the antipode of the Hopf algebras of Feynman diagrams in perturbative quantum field theory (QFT). 
The first two hold generally for arbitrary graded connected Hopf algebras and are direct consequences of the very definition of the antipode $S$ as the unique solution to $\varepsilon = S\ast I$, that is the definition of $S$ as the convolution inverse to the identity map $I$ of a Hopf algebra (here $\varepsilon$ stands for the unit of the convolution algebra of the endomorphisms of $H$).
The Dyson-Salam formula is the closed formula obtained by expanding as a formal power series in $I-\varepsilon$ the identity $I^{-1}=(\varepsilon + (I-\varepsilon))^{-1}$. 
The Bogoliubov formula is a recursive formula, obtained by rewriting the identity $\varepsilon = S\ast I$ as
$$0=S(T)+T+m\circ (S\otimes I)\otimes\DB (T),$$
where $T$ is an arbitrary element in a graded component $H_n,\ n>0$ of $H$ and $m,\DB$ stand respectively for the product and the reduced coproduct on $H$ ($\DB(T)=\Delta(T)-1\otimes T-T\otimes 1$).
The formula is
solved by induction on the degree of the graded components of $H$ and underlies the so-called Bogoliubov recursion that computes the counterterm and the renormalized values associated to the Feynman rules of a given QFT \cite{CK:2000,EMP:2000}.

The third one, the Zimmermann forest formula, is different in nature. When expanding the previous two formulas on a general Feynman diagram, terms are repeated and many cancellations occur. 
The forest formula relies on combinatorial properties that do not hold on an arbitrary graded commutative Hopf algebra, but has the nice property to reduce drastically the number of terms in the expression of the antipode; it is actually optimal in that sense.

These antipode formulas have been investigated by J.C. Figueroa and J.M. Gracia-Bondia \cite{FG1:2000,FG2:2000} in the 2000s. 
They obtained a simple direct proof of Zimmermann's formula in QFT and showed more generally that one can employ the distributive lattice of order ideals associated with a general partially
ordered set and incidence algebra techniques in order to resolve the combinatorics of overlapping divergences that motivated the development of the renormalization techniques of Bogoliubov, Dyson, Salam, Zimmermann et al.

The present article is also concerned with the forest formula, but with a different approach: we show that the formula generalizes to arbitrary right handed polynomial Hopf algebras, that is the Hopf algebras dual to the enveloping algebras of preLie algebras. This latter  structure is carried by the Hopf algebras of Feynman diagrams, but also by many other fundamental Hopf algebras since preLie algebras show up not only in QFT or related areas (statistical physics...), but also in differential geometry (an idea originating in Cayley's tree expansions), abstract algebra (Rota-Baxter algebras and operads give rise to preLie algebra structures), numerics (for the same reason as in differential geometry: differential operators give rise to preLie products), and so on. We refer e.g. to the fundational article \cite{AG:2000} and the surveys \cite{Cartier1:2000,Man:2000}.

The article is organized as follows. The first section surveys briefly enveloping algebras of preLie algebras. The second gives a direct and self-contained account of the structure of right-handed polynomial Hopf algebras. The third introduces the tree encoding of iterated coproducts and various related statistics on trees. The fourth states and exemplifies the main theorem.The last two sections are devoted to its proof.

From now on,  $k$ denotes a ground field of characteristic zero.
All the algebraic structures to be considered are defined over $k$. Since the article considers only conilpotent coalgebras, ``coalgebra'' (resp. ``Hopf algebra'') will stand for conilpotent coalgebra (resp. conilpotent Hopf algebra). 

Recall that a coalgebra $C$ with counit $\eta$ is conilpotent if and only if $\forall c\in C^+:=Ker(\eta),\exists n\in{\mathbf N}^\ast, \ \overline\Delta^{[n]}(c)=0$, where $\overline\Delta^{[n]}$ stands for the $n-1$ iterated reduced coproduct: $\overline\Delta(c):=\Delta(c)-c\otimes 1-1\otimes c$,  $\overline\Delta^{[2]}:=\overline\Delta$, and for $n\geq 3$,
$\overline\Delta^{[n]}:=(\Delta\otimes id_C^{\otimes n-2})\circ\overline\Delta^{[n-1]})$. 
In particular, cofree cocommutative coalgebra will mean cofree cocommutative in the category of conilpotent coalgebras. The cofree coalgebra over a vector space $V$ identifies then with $S(V):=\bigoplus\limits_{n\in\NM}(V^{\otimes n})^{S_n}$, the coalgebra of symmetric tensors over $V$ \cite{Haze}. Here $S_n$ stands for the symmetric group of order $n$.
It is convenient to identify elements of $S(L)$ with polynomials over $L$: we will write its elements as linear combinations of polynomials $l_1\dots l_n$ and call $S(V)$ the polynomial coalgebra over $V$. The symmetric tensor in $(V^{\otimes n})^{S_n}$ corresponding to $l_1\dots l_n$ is
$\sum\limits_{\sigma\in S_n}l_{\sigma(1)}\otimes ...\otimes l_{\sigma(n)}$.
Through this identification, the coproduct on $S(L)$ is the usual unshuffling coproduct: for arbitrary $l_1,...,l_n\in L$,
\begin{equation}\label{delta} 
\Delta(l_1...l_n)=\sum\limits_I l_I\otimes l_J
\end{equation}
where, for a subset $I$ of $[n]$, $l_I:=\prod_{i\in I}l_i$, and where $I$ runs over the (possibly empty) subsets of $[n]$ and $J:=[n]-I$.

Notice also, for further use, that, dualy, $S(V)$ is equipped with a commutative algebra structure by
the shuffle product of symmetric tensors. Through the above identification, it corresponds to the usual product of the polynomials, and we call $S(V)$ equipped with this product the polynomial algebra over $V$.

Notice finally that since conilpotent bialgebras always have an antipode (this follows e.g. from the use of the Dyson-Salam formula), the two notions of bialgebras and Hopf algebras will identify in the present article. 

\section*{Acknowledgements} The authors acknowledge support from
the grant ANR-12-BS01-0017, Combinatoire Alg{\'e}brique, R{\'e}surgence, Moules et Applications.

\section{PreLie algebras and their enveloping algebras}\label{sect:PL}

\begin{definition}
 A preLie algebra is a vector space $L$ equipped with a bilinear map $\curvearrowleft$ such that, for all $x,y,z $ in $L$:
$$(x\curvearrowleft y)\curvearrowleft z-x\curvearrowleft (y\curvearrowleft z)=
(x\curvearrowleft z)\curvearrowleft y-x\curvearrowleft (z\curvearrowleft y).$$
\end{definition}

The vector space $L$ is equipped with a Lie bracket $$[x,y]:=x\curvearrowleft y-y\curvearrowleft x,$$
see e.g. \cite{AG:2000,CL:2000} for further details.
We write $U(L)$  for the enveloping algebra of $L$ viewed as a Lie algebra. 

We will also denote by $\curvearrowleft$ the right action of the universal enveloping algebra of $L$ on $L$ that extends the pre-Lie product: $\forall a,b\in L, \ (b)a:=b\curvearrowleft a$. This action is well defined since the product $\curvearrowleft$ makes $L$ a module over $L$ viewed as a Lie algebra:
$$\forall a,b,c\in L, \ ((c)b)a-((c)a)b=(c\curvearrowleft b)\curvearrowleft a-(c\curvearrowleft a)\curvearrowleft b$$
$$=c\curvearrowleft (b\curvearrowleft a- a\curvearrowleft b)=(c)[b,a].$$

One can  equip $S(L)$, the polynomial coalgebra over $L$ with a product $\ast$ induced by $\curvearrowleft$ that makes ($S(L),\ast,\Delta)$ a Hopf algebra and the enveloping algebra of $L$ \cite{GO:2000} (recall that
an enveloping algebra carries the structure of a cocommutative Hopf algebra in which the primitive elements identify with the elements of the Lie algebra). 
The product $\ast$ is associative but not commutative and is defined as follows:
\begin{equation}\label{product}
(a_1...a_l)\ast (b_1...b_m)=\sum\limits_f B_0(a_1\curvearrowleft B_1)...(a_l\curvearrowleft B_l),
\end{equation}
where the sum is over all maps $f$ from $\{1,...,m\}$ to $\{0,...,l\}$ and $B_i:=\prod_{j\in f^{-1}(i)} b_j$. 
For example, in low degrees:
\begin{equation}\label{deg2}
a\ast b=ba+a\curvearrowleft b=ab+a\curvearrowleft b,
\end{equation}
$$a_1a_2\ast b=a_1a_2b+(a_1\curvearrowleft b)a_2+a_1(a_2\curvearrowleft b),$$
$$a\ast b_1b_2=ab_1b_2+b_1(a\curvearrowleft b_2)+b_2(a\curvearrowleft b_1)+a\curvearrowleft (b_1b_2)$$
$$=ab_1b_2+b_1(a\curvearrowleft b_2)+b_2(a\curvearrowleft b_1)+a\curvearrowleft (b_1\ast b_2 -b_1\curvearrowleft b_2)$$
$$=ab_1b_2+b_1(a\curvearrowleft b_2)+b_2(a\curvearrowleft b_1)+(a\curvearrowleft b_1)\curvearrowleft b_2 -a\curvearrowleft (b_1\curvearrowleft b_2).$$

\section{Right-handed polynomial Hopf algebras}
Notice that the increasing filtration of $S(L)$ by the degree in the previous section is respected by the product $\ast$, but the direct sum decomposition into graded components is not:
\begin{equation}\label{condit}
S^n(L)\ast S^m(L)\subset \bigoplus\limits_{n\leq i\leq n+m}S^i(L).
\end{equation}
Restricting this inclusion on the image when $n=1$, one gets a map from $L\otimes S(L)$ to $L$ which is simply the $\curvearrowleft$ map. Iterating this map and using the associativity of $\ast$ results in the identity
$$(a\curvearrowleft (b_1...b_n))\curvearrowleft (c_1...c_m)=\sum\limits_{I_1\coprod ...\coprod I_{n+1}=[m]} a\curvearrowleft ((b_1\curvearrowleft c_{I_1})...(b_n\curvearrowleft c_{I_n})c_{I_{n+1}}),$$
where the $I_i$ are possibly empty. This identity defines precisely on $L$ the structure of a symmetric brace algebra, and the two notions of symmetric brace algebras and preLie algebras happen to be equivalent (the two categories are isomorphic) \cite[Cor. 5.4]{GO:2000}.

Conversely, the categorical properties of the notion of brace algebra (see \cite{GV:2000,GJ:2000,LM:2000}) together with these results on D. Guin and M. Oudom on the symmetric brace algebra structure of the primitive part of enveloping algebras of preLie algebras \cite{GO:2000} allow to characterize enveloping algebras of preLie algebras \cite{LR:2000}. See also Turaev's \cite{turaev:2000}, where the notion of right-handed Hopf algebras used below and its links to preLie algebras and coalgebras seem to originate.

\begin{definition}
A right handed cofree cocommutative Hopf algebra is a polynomial coalgebra $S(V)$ over a vector space $V$, equipped with a product $\ast$ with unit $1\in S^0(V)=k$ that makes $S(V)$ a Hopf algebra and such that furthermore:
$$\forall n\geq 2,\ \  \ S^n(V)\ast S^m(V)\subset \bigoplus\limits_{ i\geq 2}S^i(V).$$
\end{definition}

\begin{theorem} In a right handed cofree cocommutative Hopf algebra $S(V)$,
the set of primitive elements $V$ is equipped with a preLie algebra structure by the restriction of $\ast$ to a map from $V\otimes V$ to $V$. Furthermore, $S(V)$ is the enveloping algebra of $V$ and $\ast$ identifies with the product constructed in (\ref{product}). In particular, the categories of right-handed cofree cocommutative bialgebras and the category of preLie algebras are equivalent. 
\end{theorem}

Let us give an elementary and self-contained proof of these results, disentangled of the notational complexities of brace calculus that were deviced originally for algebras up to homotopy and that are not required in this simple situation. 
We mention that the calculations in the proof are interesting on their own.

Recall first that the polynomial algebra over $V$, $S(V)=k\oplus S(V)^+$, is the free augmented commutative algebra over $V$. That is, equivalently, for an arbitrary augmented commutative algebra $A$,  $Lin(V,A^+)\cong Alg(S(V),A)$. Here, $^+$ is used to denote the augmentation ideal and $Alg$ the category of augmented commutative algebras.

By duality, for an arbitrary coaugmented cocommutative coalgebra $C=k\oplus C^+$ with coaugmentation coideal $C^+$, there is a canonical bijection (or, in categorical langage, adjunction)
$$Lin(C^+,V)\cong Coalg(C,S(V)),$$
where $Coalg$ denotes the category of coaugmented cocommutative coalgebras.
In particular, a coaugmented cocommutative coalgebra morphism to $S(V)$ is entirely characterized by its restriction to $V$ on the image.
The  inverse bijection is obtained by dualizing the isomorphism $Lin(V,A^+)\cong Alg(S(V),A)$: for $f\in Lin(C^+,V)$, the corresponding element in $Coalg(C,S(V))$ is given by 
\begin{equation}
\bigoplus\limits_{n\in\NM^\ast} f^{\otimes n}\circ \overline\Delta^{[n]},
\end{equation}
where $\overline\Delta^{[n]}$ is the iterated reduced coproduct from $C$ to $((C^+)^{\otimes n})^{S_n}\subset (C^+)^{\otimes n}$.

Let us apply this property to the Hopf algebra $S(V)$ of the Theorem.
Since $S(V)$ is a Hopf algebra, the product map $\ast$ from $C=S(V)\otimes S(V)$  to $S(V)$ is a coaugmented cocommutative coalgebra morphism, where the coproduct $\Delta$ of $S(V)\otimes S(V)$ is induced by the unshuffle coproduct(written here for notational clarity $\Delta_{S(V)}$) of $S(V)$:
\[\Delta =(I\otimes \tau\otimes I)\circ (\Delta_{S(V)} \otimes \Delta_{S(V)})\]with $\tau(x\otimes y)=y\otimes x$. 
Let us write $\pi$ for the restriction of $\ast$, on the image, to a map from $S(V)\otimes S(V)$  to $V$. By assumption, $\pi$ is null on the $V^{\otimes n}\otimes S(V),\ n\geq 2$; it is the identity map on the components $k\otimes V\cong V$ and $V\otimes k\cong V$.
For $a,b\in V$ we have $$\overline\Delta(a\otimes b)=(a\otimes 1)\otimes (1\otimes b)+(1\otimes b)\otimes (a\otimes 1)$$ and, since $\DB(a\otimes 1)=\DB(1 \otimes a)=\DB(b\otimes 1)=\DB(1 \otimes b)=0$
we get by adjunction:
\begin{equation}\label{prod2}
a\ast b=\sum_{n\geq 1} \pi^{\otimes n}\circ \DB^{[n]}(a\otimes b)=\sum_{n=1,2} \pi^{\otimes n}\circ \DB^{[n]}(a\otimes b)=\pi(a\otimes b)+a\otimes b+b\otimes a,
\end{equation}
where we recognize the equation (\ref{deg2}) expressing the associative product of two elements of a preLie algebra in the enveloping algebra in terms of the preLie product.
The same computation at higher orders would express the product $a_1\ast ...\ast a_n$ as the sum of $a_1...a_n$ and lower order terms (polynomials in $\bigoplus\limits_{k<n}S^k(V)$).

Applying this computation to the restriction to $V$ on the image of the associativity relation $\ast\circ (\ast\otimes \id)=\ast\circ (\id\otimes \ast)$, we get:
$$\forall a,b,c\in V,\ \pi (\pi(a\otimes b)\otimes c)=\pi(
a \otimes (\pi(b\otimes c)+b\otimes c+c\otimes b),$$
or 
\begin{equation}\label{symmetry}
\pi (\pi(a\otimes b)\otimes c)-\pi(
a \otimes (\pi(b\otimes c))=\pi (a\otimes (b\otimes c+c\otimes b)).
\end{equation}
Since the last expression is symmetric in $b$ and $c$, we get finally
$$\pi (\pi(a\otimes b)\otimes c)-\pi(
a \otimes (\pi(b\otimes c))=\pi (\pi(a\otimes c)\otimes b)-\pi(
a \otimes (\pi(c\otimes b)),$$
where we recognize the preLie identity.

The same calculation can be repeated with higher tensor products: the restriction  to $V$ on the image of the product map from $S(V)^{\otimes n}$ to $S(V)$ can be computed on a tensor product of elements of $V$ as
$$\pi (...\pi(\pi(a_1\otimes a_2)\otimes a_3)...\otimes a_n),$$
or as $$\pi (a_1\otimes (a_2\ast ...\ast a_n))=\pi (a_1\otimes ((a_2 ... a_n)+ \ R),$$
where the remainder term $R$ is a polynomial in $\bigoplus\limits_{k<n}S^k(V)$.
This shows that the restrictions of $\pi$ to $V\otimes S^n(V)$ can be computed inductively and depend only on the value of $\pi$ on $V\otimes V$.

The  Theorem follows: $V$, the Lie algebra of primitive elements of $S(V)$ is a preLie algebra. By the Cartier-Milnor-Moore theorem \cite{Cartier2:2000,mm:2000,Patras1:2000}, $S(V)$ is its enveloping algebra (the proof given in \cite{Patras:2000,Patras1:2000} requires only cocommutativity and conilpotency, not the graduation hypothesis of \cite{mm:2000}). Since the preLie algebra structure of $V$ determines uniquely the algebra structure of $S(V)$, the product $\ast$ identifies with the product computed using (\ref{product}).

The previous results dualize (\cite{turaev:2000}). 
\begin{definition}
A right handed polynomial Hopf algebra is a polynomial algebra $S(V)$ over a vector space $V$, equipped with a coproduct $\Delta$ (with counit the projection from $S(V)$ to $S^0(V)=k$) that makes $S(V)$ a Hopf algebra and such that furthermore the coproduct $\Delta$ is right-handed, that is $$\overline\Delta(V)\subset V\otimes S(V).$$
\end{definition}

\begin{proposition} Let $S(V)$ be a right-handed polynomial Hopf algebra.
Then, $V$ is equipped with a preLie coalgebra structure by the restriction of $\Delta$ to a map from  $V$ to  $V\otimes V$. 
\end{proposition}

The structure theorems for right handed cofree cocommutative Hopf algebras dualize perfectly in the framework of Milnor-Moore \cite{mm:2000}, that  is when the $S(V)$ are furthermore connected graded Hopf algebra (in general, dualizing structure theorems for Hopf algebras requires some care since the dual of a coalgebra in an algebra, but the converse is not true in general -the notion of restricted duals has to be used. In the category of graded vector spaces, these difficulties disappear).
By graded duality, graded connected right-handed polynomial bialgebras are coenveloping coalgebras of graded connected preLie coalgebras and the corresponding categories are equivalent.

Most combinatorial Hopf algebras are graded (see \cite{Cartier2:2000}), and such are the ones of Feynman diagrams in QFT, where the diagrams are graded by their number of loops \cite{CK:2000}.

\section{Toward chains and forests}\label{sec:chtex}
Chains and forests are naturally associated to the action of the iterated reduced coproduct on the Feynman graphs of a given QFT: chains are successions of strict inclusions of subgraphs, whereas forests are families of subgraphs satisfying certain technical conditions, essentially such that the connected components of the subgraphs in a chain form a forest.
We refer to \cite{FG2:2000} for a detailed analysis of these notions, that have also appeared recently in relation to Hopf algebra structures in control theory \cite{DEG:2000}. 

This section aims at defining the analog notion of a forest (indeed trees) in the more general context of right-handed polynomial bialgebras. As we shall see, some tree indexations naturally appear in the computation of iterated coproducts.

Let from now on in this article $H=S(V)$ be a conilpotent right-handed polynomial Hopf algebra. We assume that $V$ has a basis $\mathcal B=\{ b_i\}_{i\in X}$, where $X=\{1,\dots,n\}$ or $X=\NM^\ast$, fixed once for all. Notice that in most application domains of the theory of preLie algebras and coalgebras (see \cite{Cartier1:2000,Man:2000}), there is a natural choice for the basis $\mathcal B$ and therefore a natural notion of chain and forest will result.

In the sequel, we consider multisets over $X$ and, for any such multisets $I,J$, write $I\cup J$ for their union. For example, $\{1,2,2\}\cup\{2,3,3\}=\{1,2,2,2,3,3\}$.
With these notations, one can consider monomials $b_I=\prod_{i\in I}b_i$, so that $b_I.b_J=b_{I\cup J}$ and one can note $b_\emptyset$ the unit of $H$.
Let us fix $b_i\in\mathcal B$. We aim at computing the value of the antipode on $b_i$, $S(b_i)$, that can already be expressed using the Dyson-Salam formula with the help of iterated coproducts:
\[ S(b_i)=\sum_{k\geq 1} (-1)^k m^{[k]}\circ \DB^{[k]}(b_i),\]
where $m^{[k]}$ is the $k$--fold product.

We can first expand the reduced coproduct of $b_i$ as follows:
\[\DB (b_i)=\sum_{i_0, I\not=\emptyset} \lambda^{i ; i_0}_{I} b_{i_0}\otimes b_{I}.
\]
The coefficients $\lambda^{i_0 ; i}_{I}$ completely determine the  coproduct and its action on products, as well as the action of the iterated coproducts. At this stage, one can opt for a graphical representation and consider that, in $\DB (b_i)$, the different terms are indexed by non planar decorated corollas whose root is decorated by $(i ;i_0)$ and leaves decorated by the positive integers:
\[
\DB (b_i)=\sum \lambda\Big(\corollaa \Big) b_{i_0}\otimes b_{i_1}\dots b_{i_k}.
\]
Here, non planar means as usual for trees that the ordering of the branches does not matter, reflecting the commutativity of the product.

Let us now consider a single term in this reduced coproduct, for instance
\[
\lambda\Big (\corollab \Big ) b_{i_0}\otimes b_{i_1} b_{i_2}.
\]
One can observe that 
\[
\DB(b_{i_1} b_{i_2}) = (1\otimes b_{i_1}+b_{i_1} \otimes 1+\DB(b_{i_1}))(1\otimes b_{i_2} +b_{i_2}\otimes 1+\DB(b_{i_2}))-1\otimes b_{i_1}b_{i_2}-b_{i_1}b_{i_2}\otimes 1
\]
so that the contribution of $\lambda\left(\corollab \right) b_{i_0}\otimes b_{i_1} b_{i_2}$ to $\DB^{[3]}(b_i)=(Id\otimes \DB)(\DB(b_i))$ will split in four terms, whose complexity is encoded by the appearance of products of coefficients $\lambda^{\cdot;\cdot}_{\cdots}$.  

There is a first term with no more "complexity" than in $\DB(b_i)$:
\[\lambda
\Big (\corollab \Big ) b_{i_0}\otimes (b_{i_1}\otimes b_{i_2} +b_{i_2}\otimes  b_{i_1}).
\]
There is a second term, where only the reduced coproduct of $b_{i_1}$ occurs:
\[\lambda\Big (\corollab \Big ) (\sum \lambda^{i_1;i_{1,0}}_{i_{1,1}\dots i_{1,k}} b_{i_0}\otimes (b_{i_{1,0}}\otimes b_{i_2}b_{i_{1,1}\dots i_{1,k}} +b_{i_2}b_{i_{1,0}}\otimes  b_{i_{1,1}\dots i_{1,k}})),
\] and this contribution is naturally indexed by the trees:
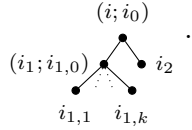
\[\corollac .\]
In the same way, there is a contribution, corresponding to $(1\otimes b_{i_1}+b_{i_1} \otimes 1)\DB(b_{i_2})$ indexed by the trees:
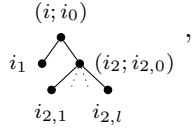
\[\corollad ,\] and finally the terms in relation with $\DB(b_{i_1})\DB(b_{i_2})$ that will be indexed by trees 
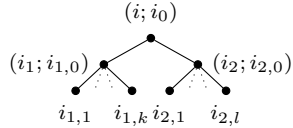
\[\corollae .\]

When iterating the reduced coproduct, such groups of terms naturally appear, labeled by trees that encode the presence of the coefficients $\lambda^{\cdot;\cdot}_{\cdots}$.

\begin{definition}\label{def:tree}
Let us consider   finite rooted trees (connected and simply connected finite posets with a unique minimal element) whose internal vertices are decorated by pairs $p=(p_1;p_2)$ of positive integers and leaves are decorated by positive integers (note that in the tree with only one vertex, the vertex is considered as a leaf).

 A forest is simply a commutative product of such trees. 
For any internal vertex $x$ ($x\in Int(T)$) in such tree or forest, we note $d(x)=(d_1(x);d_2(x))$ its decoration, and, if $x$ is a leaf ($x\in Leaf(T)$), we note for convenience its decoration $d(x)=d_1(x)=d_2(x)$. For any internal vertex x, we also note, $succ(x)$ the set of its immediate successors.

If the root of a tree is decorated by $i$ or $(i;i_0)$, we say that the tree is associated to $b_i$ ($T\in \mathcal{T}_{i}$).
For a given pair $p$, we note $B^+_p(T_1\dots T_s)$ the tree obtained by adding a common root decorated by $p$ to the trees $T_1\dots T_s$.
\end{definition}

The reader will notice that in the QFT terminology such trees are called forests (because to each tree is associated a forest by cutting the root).

\begin{definition}
The length of a tree, $l(T)$, is the number of elements in $T$ viewed as a poset. The height $h(T)$ of a tree is the maximum number of elements in a chain from the root to a leaf. 

The coefficient $\lambda(T)$ is defined as follows: $\lambda(\bullet_i)=1$ and if $T=B^+_{(i;i_0)} (T_1\dots T_s)$, then 
\[\lambda(T)=\lambda^{i,i_0}_{i_1,\dots,i_s}\lambda(T_1)\dots\lambda(T_s)\]when $T_1,\dots, T_s$ are respectively associated to $b_{i_1},\dots,b_{i_s}$. In other words
\[ \lambda(T)=\prod_{x\in Int(T)}\lambda^{d(x)}_{d_1(succ(x))}.\]

The $b$-value of a tree is the element $v(T)$ of $S(V)$ defined as 
\[v(T)=\prod_{x\in T} b_{d_2(x)}\]
\end{definition} 

We extend naturally these notions to forests:
\begin{eqnarray*}
l(T_1\dots T_s) &=& l(T_1)+\dots +l(T_s)\\
h(T_1\dots T_s) &=& max(h(T_1),\dots h(T_s)) \\
\lambda(T_1\dots T_s)&=&\lambda(T_1)\dots \lambda(T_s) \\
v(T_1\dots T_s)&=&v(T_1)\dots v(T_s).
\end{eqnarray*}

In the sequel, once a forest $F$ is given, we will note abusively, for any vertex $x$ of $F$, $b_x=b_{d_2(x)}$. 

For example, for the tree $$T=\exarb ,$$ $l(T)=8$, $h(T)=3$, 
$ \lambda(T)=\lambda^{i;i_{0}}_{i_{1},i_{2}}
\lambda^{i_{1};i_{1,0}}_{i_{1,1},i_{1,2}}
\lambda^{i_{2};i_{2,0}}_{i_{2,1},i_{2,2},i_{2,3}}$ and $$v(T)=b_{i_{0}}b_{i_{1,0}}b_{i_{1,1}}b_{i_{1,2}}
b_{i_{2,0}}b_{i_{2,1}}b_{i_{2,2}}b_{i_{2,3}}.$$

We can also rephrase the definition of $\DB$:
\begin{equation}\label{cop:deux}
\DB(b_i)=\sum_{T=B^+_{(i;i_0)}(\bullet_{i_1}\dots\bullet_{i_k}) \in \mathcal{T}_i ; h(T)=2} \lambda(T) v(\bullet_{i_0})\otimes v(\bullet_{i_1}\dots \bullet_{i_k}).
\end{equation}
We can now state the Zimmermann forest formula in the framework of right-handed polynomial algebras.

\section{The PreLie forest formula}
\begin{theorem}\label{antipode}
The value of the antipode $S$ of the right-handed polynomial bialgebra $S(V)$ on an element $b_i\in \mathcal{B}$ is given by the cancellation free formula:
\begin{equation}
S(b_i)=\sum\limits_{T\in \mathcal{T}_i}(-1)^{l(T)}\lambda(T)v(T)
\end{equation}
\end{theorem}

By cancellation free, we refer to the fact that each tree appears only once, as in the classical QFT Zimmermann's forest formula.
Several terms corresponding to the same tree would instead appear in the Dyson-Salam (and Bogoliubov) formula, as illustrated below.

We postpone the proof to the next section. Let us first show on an elementary example how the notion of forest and the forest formula behave concretely. We consider the emblematic case (see e.g. \cite{FG2:2000}) of the Fa\`a di Bruno Hopf algebra encoding the substitution product in the algebra of
 formal power series $$f(t) =t+
\sum\limits_{n=2}^\infty f_n \frac{t^n}{n!}.$$
On the polynomial algebra generated by the coordinate functions $a_{n}(f) := f_n , n \geq 2$, the substitution product translates into the coproduct
$$\Delta(a_n)=\sum\limits_{k=1}^na_k\otimes B_{n,k}(a_1,...,a_{n+1-k}),$$
where $a_1:=1$ and the $B_{n,k}$ are the (partial, exponential), Bell polynomials defined by the series expansion
$$\exp (u \sum\limits_{m\geq 1} x_m \frac{t^m}{m!})= 1+\sum\limits_{n\geq 1}
\frac{t^n}{n!}[\sum\limits_{k=1}^n u^kB_{n,k}(x_1, . . . , x_{n+1-k})].$$
Setting $b_n:=a_{n+1}$, we get
$$\overline\Delta(b_1)=0,\ \overline\Delta(b_2)=3b_1\otimes b_1,\ \overline\Delta(b_3)=6b_2\otimes b_1+b_1\otimes (3b_1^2+4b_2).$$

To compute $S(b_3)$, let us apply the classical Dyson-Salam formula obtained by expanding as a formal power series in $I-\varepsilon$ the identity $S=I^{-1}=(\varepsilon + (I-\varepsilon))^{-1}=\sum_n(-1)^n (I-\varepsilon)^{\ast n}$. 
We get, grouping the terms according to the powers $(I-\varepsilon)^{\ast n}$ and (since our goal is here to understand the structure of the calculation of the antipode on an example) avoiding to identify immediately the terms inside these groups:
$$S(b_3)=-b_3+(6b_2\cdot b_1+3\ b_1\cdot b_1^2+4\ b_1\cdot b_2)-(3\cdot 2\ b_1\cdot b_1^2+4\ \cdot 3\  b_1\cdot b_1^2)$$
$$=-b_3+10b_1b_2-15b_1^3.$$
Overall, 6 terms appear that can be resummed into three terms.

On the other hand
$\overline\Delta(b_3)=6b_2\otimes b_1+b_1\otimes (3b_1^2+4b_2)$
can be rewritten
$$\overline\Delta(b_3)=\lambda(\tfdbA) b_2\otimes b_1+\lambda(\tfdbB)
b_{1}\otimes b_1 b_1
+\lambda(\tfdbC)
b_1\otimes b_2.$$
When iterating the coproduct, as in the previous section, we get:
\[
\DB^{[3]}(b_3)= 2 \lambda(\tfdbB)
b_{1}\otimes b_1 \otimes b_1
+\lambda(\tfdbC)\lambda(\tfdbD)
b_1\otimes b_1 \otimes b_1
\]
where 
\[
\lambda(\tfdbC)\lambda(\tfdbD)=4.3=12=\lambda(\tfdbE)
\]

Overall, we have now 5 forests: $\bullet_3$, $\tfdbA$, $\tfdbB$, $\tfdbC$, $\tfdbE$, and instead of the six terms in the formula $S(b_3)=-b_3+m^{[2]}\circ \DB^{[2]}(b_3)-m^{[3]}\circ \DB^{[3]}(b_3)$, the indexation by trees, that corresponds to the forest formula gives
\[
S(b_3)= -b_3 + \lambda(\tfdbA) b_2\cdot b_1- \lambda(\tfdbB)
b_{1}\cdot b_1 \cdot b_1
+\lambda(\tfdbC)
b_1\cdot b_2- \lambda(\tfdbE)b_{1}\cdot b_1 \cdot b_1
\]
That contains only 5 terms and takes into account cancellations between the terms associated to $\tfdbB$.

Although extremely elementary, this example gives the flavour of the general pattern followed by the forest formula and of the cancellations occurring. The reader can also compare with other examples as in \cite{FG2:2000} or \cite{DEG:2000}.

\section{Chains and linearization of forests}

This section aims at proving the main theorem, Thm \ref{antipode} and, in the process, introduces various useful tools in order to understand the behaviour of right handed polynomial bialgebras $S(V)$.

The first computations we did in section \ref{sec:chtex} suggest that, when iterating the reduced coproduct, the tensor products  we get can be associated to trees, so that, in the end, some cancellations occur in the computation of the antipode and yields the preLie forest formula.

The following notions, inspired by analogous constructions on finite topologies (a generalization of posets) \cite{FMP:2000}, aim at encoding these formulas.

\begin{definition}
Let $P$ be a finite poset of cardinality $n$. A $k$-linearization of $P$ is a surjective, order preserving map $f$ from $P$ to $[k]$, where $k\leq n$. We write $f\in k-lin (\mathcal P)$.
\end{definition}

Here, order preserving means that strict inequalities are preserved: $x<y$ implies $f(x)<f(y)$. Note that $k$ must be greater or equal to the "height" of the poset, that is the length of its maximal interval.

\begin{definition}
Let $P$ be a forest (as in Definition \ref{def:tree}) with a given decoration $d=(d_1,d_2)$. If $f$ is a $k$-linearization of $P$, we call the tensor product
$$C(f):=(\prod\limits_{x_1\in f^{-1}(1)}b_{x_1})\otimes ...\otimes (\prod\limits_{x_k\in f^{-1}(k)}b_{x_k})$$
a $k$-chain of $P$. As before, $b_{x_i}$ stands for $b_{d_2(x_i)}$.

\end{definition}

Since, by their very definition, the decorations of trees associated to $b_i$  run over all the indices of basis elements appearing in the various iterated reduced coproducts of $b_i$ (ordered according to their relative positions in the iterations of the reduced coproducts), a fundamental key to the proof is that $k$-linearizations  describe all the tensors of length $k$ that can be obtained in the $k$-fold iterated reduced coproduct.

The proof of theorem will follow from the two following fundamental lemmas, whose proof is postponed to the next Section:
 
\begin{lemma}\label{lem:prodcop}
 Let $I=\{i_1,\dots,i_s\}$ be a multiset and $\mathcal{F}_I=\{T_1\dots T_s, T_j\in \mathcal{T}_{i_j}\}$, then
 \[
 \DB(b_I)=\sum_{F\in \mathcal{F}_I} \sum_{f\in 2-lin(F)}\lambda(F)C(f)
 \]
 
\end{lemma}

\begin{lemma} \label{lem:itcop}
We have, for the action of the $k$-fold iterated coproduct $\overline\Delta^{[k]}:=(\id^{\otimes k-2}\otimes \overline\Delta)\circ ...\circ (\id \otimes \overline\Delta)\circ\overline\Delta$:
$$\overline\Delta^{[k]}(b_i)=\sum\limits_{T\in \mathcal{T}_i}\sum\limits_{f\in k-lin (T)}\lambda(T)C(f).$$

\end{lemma}

As a corollary, 
\begin{corollary}
For $b_i\in \mathcal B$, we have:
$$(\id-\varepsilon)^{\ast k}(b_i)=\sum\limits_{T\in \mathcal{T}_i}\sum\limits_{f\in k-lin (\mathcal F)}\lambda(T)v(T),$$
and
$$S(b_i)=\sum\limits_{k\geq 1}
\sum\limits_{T\in \mathcal{T}_i}\sum\limits_{f\in k-lin (T)}(-1)^k\lambda(T)v(T).$$
\end{corollary}

The proof of Thm \ref{antipode} boils down in the end therefore to proving that, for a given tree $T$, 
$$\sum\limits_{k\geq 1}\sum\limits_{f\in k-lin (T)}(-1)^k =(-1)^{l(T)}.$$
The identity follows from the general Proposition:

\begin{proposition}
For an arbitrary rooted tree $P$ of cardinality $m$ viewed as a poset, we have:
$$\sum\limits_{f\in k-lin (P)}(-1)^k =(-1)^m.$$
\end{proposition}

Let us prove the Proposition by induction on $m$: it is obvious when $m=1$. Let us assume that the Proposition is true for posets of cardinality less or equal $m$. A poset $P'$ of cardinality $m+1$ can always be written $P'=P\cup\{x\}$, where $x$ is a maximal element in $P'$. Let us also introduce the predecessor $y$ of $x$ for the tree structure (the maximal $z$ with $z<x$).

A linearization of $P'$ can be obtained from a $k$-linearization $f$ of $P$ as follows (all linearizations of $P'$ are obtained in that way).

Consider the sequence $(F_1,\dots,F_k)= (f^{-1}(1),\dots, f^{-1}(k))$ with $p=f(y)$ ($y\in F_p$). The linearizations that can be obtained by inserting $x>y$  correspond to the sequences:
\[
\begin{array}{ll}
 (F_1,\dots,F_p,\{x\},F_{p+1},\dots ,F_k) & (F_1,\dots,F_p,F_{p+1}\cup \{x\}, F_{p+2},\dots ,F_k) \\
 (F_1,\dots,F_p,F_{p+1},\{x\},F_{p+2},\dots ,F_k) & (F_1,\dots,F_p,F_{p+1}, F_{p+2}\cup \{x\},\dots ,F_k) \\
 \vdots &\vdots \\
 (F_1,\dots ,\{x\}, F_k) &(F_1,\dots ,F_k\cup \{x\}) \\
 (F_1,\dots ,F_k,\{x\})
 \end{array}
 \]
%
We get that $f$ gives rise to $(k-p)$ $k$-linearizations of $P'$ and $(k-p+1)$ $k+1$-linearizations of $P'$. 
Finally:
$$\sum\limits_{f\in k-lin (P')}(-1)^k =\sum\limits_{f\in k-lin (P)}(-1)^k((k-p)-(k-p+1))$$
$$=-\sum\limits_{f\in k-lin (P)}(-1)^k=
(-1)^{m+1}.$$

\section{Iterated coproducts and trees}

We postponed to this section the proof of lemmas \ref{lem:prodcop} and \ref{lem:itcop}, the first one serving in the proof of the latter.

\subsection{Proof of Lemma \ref{lem:prodcop}}
We want to prove that for a given multiset $I=\{i_1,\dots,i_s\}$ and $\mathcal{F}_I=\{T_1\dots T_s, T_j\in \mathcal{T}_{i_j}\}$, we have
 \[
 \DB(b_I)=\sum_{F\in \mathcal{F}_I} \sum_{f\in 2-lin(F)}\lambda(F)C(f)
 \]
 
First observe that, as linearizations are strictly increasing, $2-lin(F)$ is empty if $h(F)>2$. 

The proof is recursive on the cardinal of $I$. 

If $I=\{i\}$, we recover the formula (\ref{cop:deux}): Let $T\in \mathcal{T}_i$, if $h(T)=1$, then $T=\bullet_i$ and $2-lin(\bullet_i)$ is empty, otherwise, $h(T)=2$ and $T=B^+_{(i_0;i)}(\bullet_{i_1}\dots\bullet_{i_k})$. In this latter case, the only possible $2$-linearization sends the root on 1 and all the other vertices on 2.

Suppose now that the result holds for any monomial $b_I=b_{\{i_1,\dots, i_s\}}$ of a given length $s$. For the monomial $b_I.b_j$, the right-hand term $S$ of the above formula can be written
\[
S=\sum_{F\in \mathcal{F}_I}\sum_{T\in \mathcal{T}_j} \sum_{f\in 2-lin(F.T)}\lambda(F.T)C(f)
\]
and $\lambda(F.T)=\lambda(F)\lambda(T)$. For the poset $F.T$, no vertex of $F$ is comparable to a vertex of $T$ and any $2$-linearization $f$ of $F.T$ can be deduced from the restrictions of $f$ to $F$ and $T$. For $k=1,2$ let $F^k=f^{-1}(k)\cap F$ and $T^k=f^{-1}(k)\cap T$, so that $f^{-1}(k)$ is the disjoint union of $F^k$ and $T^k$:
\begin{enumerate}
\item If $f$ is such that none of the $F^k$ or $T^k$ is empty, then the restrictions of $f$ to $F$ and $T$ are 2-linearizations. The set of such $f$ is in bijection with $2-lin(F)\times 2-lin(T)$ and the corresponding subsum of S is:
\[S_1= \sum_{F\in \mathcal{F}_I \atop T\in \mathcal{T}_j} \sum_{f_1\in 2-lin(F)\atop
f_2\in 2-lin(T)}\lambda(F)\lambda(T)C(f_1)C(f_2)= \DB(b_I) \DB(b_j)\]
\item If $f$ is such that $F^1$ and $F^2$ are nonempty but $(T^1,T^2)= (T^1,\emptyset)$, necessarily, no elements of $T^1$ are comparable: otherwise two such elements could not be in the same subset $f^{-1}(1)$. We thus have $T^1={\bullet_j}$ and, finally, the corresponding subsum is:\[
S_2=\sum_{F\in \mathcal{F}_I} \sum_{f_1\in 2-lin(F)}\lambda(F) C(f_1).(b_j\otimes 1)=\DB(b_I).(b_j\otimes 1)\]
\item For the same reason, if $f$ is such that $F^1$ and $F^2$ are nonempty but $(T^1,T^2)= (\emptyset ,T^2)$, necessarily, $T^2={\bullet_j}$ and the corresponding subsum is:\[
S_3=\sum_{F\in \mathcal{F}_I} \sum_{f_1\in 2-lin(F)}\lambda(F) C(f_1).(1\otimes b_j)=\DB(b_I).( 1 \otimes b_j )\]
\item If $f$ is such that $T^1$ and $T^2$ are nonempty but $(F^1,F^2)= (F^1,\emptyset)$, necessarily, $F^1=\bullet_{i_1}\dots \bullet_{i_s}$ and the corresponding subsum is:\[
S_4=\sum_{T\in \mathcal{T}_j} \sum_{f_2\in 2-lin(T)}\lambda(T) C(f_2).(b_J\otimes 1)=\DB(b_j).(b_I\otimes 1)\]
\item If $f$ is such that $T^1$ and $T^2$ are nonempty but $(F^1,F^2)= (\emptyset ,F^2)$, necessarily, $F^2=\bullet_{i_1}\dots \bullet_{i_s}$ and the corresponding subsum is:\[
S_5=\sum_{T\in \mathcal{T}_j} \sum_{f_2\in 2-lin(T)}\lambda(T) C(f_2).(1\otimes b_I)=\DB(b_j).( 1 \otimes b_I )\]
\item If $(T^1,T^2)= (T^1,\emptyset)$ and $(F^1,F^2)= (\emptyset ,F^2)$, this correspond to a unique $2$-linearization that gives : \[
S_6=b_j\otimes b_I\]
\item If $(T^1,T^2)= (\emptyset ,T^2)$ and $(F^1,F^2)= (F^1,\emptyset)$, this correspond to a unique $2$-linearization that gives : \[
S_7=b_I\otimes b_j\]
\end{enumerate}
Putting together these seven sums, this gives\[
S+1\otimes b_{I\cup j}+ b_{I\cup j}\otimes 1=(\DB(b_I)+1 \otimes b_I+  b_I\otimes 1)(\DB(b_j)+1 \otimes b_j+  b_j\otimes 1)
 \]
Thus \[
S=\Delta(b_{I\cup j})-1\otimes b_{I\cup j}- b_{I\cup j}\otimes 1=\DB(b_{I\cup j})
\] This ends the proof.

\subsection{Proof of Lemma \ref{lem:itcop}}

It remains to prove that 
$$\overline\Delta^{[k]}(b_i)=\sum\limits_{T\in \mathcal{T}_i}\sum\limits_{f\in k-lin (T)}\lambda(T)C(f).$$

We already proved the above formula for $k=2$ (see formula \ref{cop:deux}). Suppose the result holds for a given $k\geq 2$ and consider the right-hand side of the formula at order $k+1$:
$$R= \sum\limits_{T\in \mathcal{T}_i}\sum\limits_{f\in (k+1)-lin (T)}\lambda(T)C(f).$$
Note that we can restrict the sum to trees of height greater than $k$. 

\begin{definition}
A nonempty corolla cut $C$ of a tree $T$ ($C\in Ccut(T)$) is a subset of $T$ such that: (1) its elements are maximal or predecessors of maximal elements, (2) if $y\in C$ then $\{x\ ;\ x>y\}\subset C$.
\end{definition}
 Such a cut inherits the decorations of $T$ and the order induced by $T$. It is clear that this is a forest of height 1 or 2. 
The reader can easily check that for any corolla cut $C$ of a tree $T$, $\lambda(T)=\lambda(T/C)\lambda(C)$. 
 
 For instance, if 
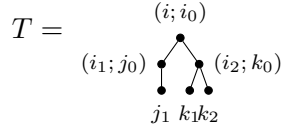
\[T=\begin{arb}
\rd{(i;i_0)}
child{\vl{(i_1 ; j_{0})} child{\vb{j_{1}}} } child{\vr{(i_2 ; k_{0})}  child{\vb{k_{1}}} child{\vb{k_{2}}} } ;
\end{arb}
\]
We obtain 7 corolla cuts of height 1 by choosing  1,2 or 3 leaves. As for the corolla cuts of heights 2, we get:
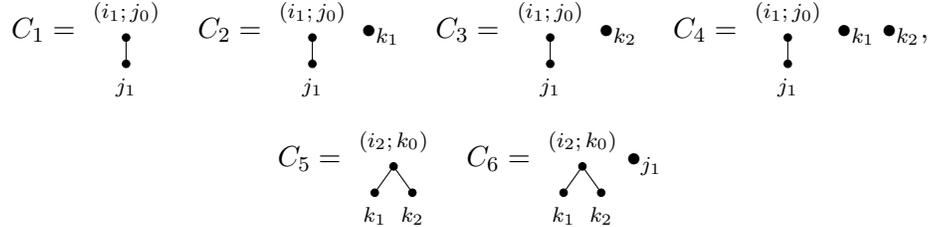
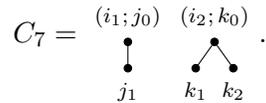
\[
C_1= \begin{arb}
\rd{(i_1 ; j_{0})}
child{\vb{j_{1}}};
\end{arb}
\quad C_2=\begin{arb}
\rd{(i_1 ; j_{0})}
child{\vb{j_{1}}};
\end{arb} 
\bullet_{k_1}
\quad C_3=\begin{arb}
\rd{(i_1 ; j_{0})}
child{\vb{j_{1}}};
\end{arb} 
\bullet_{k_2}
\quad C_4=\begin{arb}
\rd{(i_1 ; j_{0})}
child{\vb{j_{1}}};
\end{arb} 
\bullet_{k_1}  \bullet_{k_2},
\]
\[
C_5 =\begin{arb}
\rd{(i_2 ; k_{0})}
child{\vb{k_{1}}} child{\vb{k_{2}}}  ;
\end{arb}
\quad C_6=\begin{arb}
\rd{(i_2 ; k_{0})}
child{\vb{k_{1}}} child{\vb{k_{2}}}  ;
\end{arb}
\bullet_{j_1}
\]
and, finally,
\[
C_7=\begin{arb}
\rd{(i_1 ; j_{0})}
child{\vb{j_{1}}};
\end{arb} 
\begin{arb}
\rd{(i_2 ; k_{0})}
child{\vb{k_{1}}} child{\vb{k_{2}}}  ;
\end{arb}.
\]
Once such a corolla cut $C=T_1..T_s$ is given, we note $T/C$ the tree obtained as follows: for $1\leq i\leq k$, if $h(T_i)= 2$ (a "true terminal" corolla) remove all the maximal elements of $T_i$ in $T$ and replace the decoration $(k;l)$ of the root of $T_i$  by $k$ in the new tree. In the previous example $T/C=T$ if $h(C)=1$,
\[T/C_1=T/C_2=T/C_3=T/C_4=\begin{arb}
\rd{(i;i_0)}
child{\vl{i_1 } } child{\vr{(i_2 ; k_{0})}  child{\vb{k_{1}}} child{\vb{k_{2}}} } ;
\end{arb},
\]
\[T/C_5=T/C_6=\begin{arb}
\rd{(i;i_0)}
child{\vl{(i_1 ; j_{0})} child{\vb{j_{1}}} } child{\vr{i_2 }  } ;
\end{arb}
\]and, finally
\[
T/C_7=\begin{arb}
\rd{(i;i_0)}
child{\vl{i_1 }  } child{\vr{i_2 }  } ;
\end{arb}.
\]
We write $T\wedge C$ for the set of leaves of $T/C$ that coincides with minimal elements of $C$ and
for a given $T$ and $g\in m_{lin}(T)$, write 
$$(T^1_g,\dots,T^m_g)=(g^{-1}(1),\dots,g^{-1}(m)),$$ 
and $$C(g)=C_1(g)\otimes \dots \otimes C_m(g).$$

Let us observe that, for any $k+1$-linearization $f$ of a tree $T$, since $f$ is strictly increasing, 
\begin{itemize} 
\item The set $C=T^k_f \cup T^{k+1}_f$ is a nonempty corolla cut of $T$.
\item The restriction of $f$ to this cut determines a unique $2$-linearization $f_C$ of $C$.
\item The map $f^C$ defined on $T/C$ by $f^C(x)=f(x)$ if $f(x)<k$ and $f^C(x)=k$ otherwise is a $k$ linearization of $T/C$ and $(f^C)^{-1}(k)=T\wedge C$.
\end{itemize}

Conversely, on can  associate to a sequence $(C,g,h)\in Ccut(T)\times k-lin(T/C)\times 2-lin(C)$ such that $g^{-1}(k)=T\wedge C$  a unique $(k+1)$-linearization   on $T$ given by the ordered partition 

$$(g^ {-1}(1), \dots, g^{-1}(k-1),h^{-1}(1),h^{-1}(2)).$$

Thanks to this bijection, $R$ is equal to
$$
\sum_{\begin{array}{c}
T\in \mathcal{T}_i \\
C\in Ccut(T) \\ 
g\in k-lin(T/C) ; g^{-1}(k)=T\wedge C  \\ 
h\in 2-lin(C)
\end{array} } \lambda(T/C)\lambda(C)C_1(g)\otimes\dots\otimes C_{k-1}(g)\otimes C_1(h)\otimes C_2(h)
$$

We can reindex this sum by $T'=T/C$ that run over $\mathcal{T}_i$, $g\in k-lin(T')$, $C$ is a forest of height lower or equal to 2 whose set of roots $\{r_1,..., r_s\}$ has the same cardinal than $g^{-1}(k)$ (namely the previous cardinality of $T\wedge C$) and their decoration $d_1$ coincide. If we note $I=d_1(g^{-1}(k))$ :
$$
R=\sum_{ \begin{array}{c}
T'\in \mathcal{T}_i \\
g\in k-lin(T)  \\ 
C\in \mathcal{F}_I \\ 
h\in 2-lin(C)
\end{array} } \lambda(T')\lambda(C)C_1(g)\otimes \dots\otimes C_{k-1}(g)\otimes C_1(h)\otimes C_2(h)
$$
and, since $C_k(g)=b_I$, we get, thanks to the previous lemma:
$$
R=\sum_{\begin{array}{c}
T'\in \mathcal{T}_i \\

g\in k-lin(T)  \\ 

\end{array} } \lambda(T')C_1(g)\otimes\dots\otimes C_{k-1}(g)\otimes \DB(C_{k}(g))=(Id^{\otimes^{k-1}}\otimes \DB)\circ \DB^{[k]}(b_i).
$$

This ends the proof of the lemma.


\end{document}